\documentstyle[11pt]{article}
\begin{document}
 \baselineskip 0.22in

\title{\textbf{Equivalence and Strong Equivalence between Sparsest and Least $\ell_1$-Norm Nonnegative
  Solutions of Linear Systems  and Their Application }}
\author{YUN-BIN ZHAO \thanks{School of Mathematics, University of Birmingham,
Edgbaston, Birmingham B15 2TT,  United Kingdom ({\tt
y.zhao.2@bham.ac.uk}).  This work was supported by the Engineering
and Physical Sciences Research Council (EPSRC) under the grant
\#EP/K00946X/1. }}

  \date{{\small (First version: June 2012; Revised in November 2012, September and December 2013)}}

\maketitle

{\bf Abstract.}  Many practical problems can be formulated as
  $\ell_0$-minimization problems with nonnegativity constraints, which
 seek the sparsest nonnegative solutions to underdetermined  linear
systems. Recent study indicates that
    $\ell_1$-minimization is   efficient    for solving some classes of
   $\ell_0$-minimization problems.
 From a mathematical point of view, however,
the understanding of the relationship between $\ell_0$- and
 $\ell_1$-minimization  remains incomplete. In this paper, we
further
 discuss several theoretical questions associated with these two problems.  For instance,
  how to completely characterize the uniqueness of least
  $\ell_1$-norm nonnegative solutions to a linear system, and  is there any
  alternative
   matrix property that is different
  from existing ones, and
  can fully characterize the uniform recovery of $K$-sparse nonnegative vectors?
  We prove that the fundamental strict complementarity theorem of
linear programming can  yield a necessary and sufficient condition
for a linear system  to have a unique least $\ell_1$-norm
nonnegative solution. This condition leads naturally to the
so-called range space property (RSP)  and the `full-column-rank'
property, which altogether provide a   broad understanding of
the relationship between $ \ell_0$- and
 $\ell_1$-minimization. Motivated by these results, we introduce the concept of the `RSP of
 order $K$' that turns out to be a full characterization of the  uniform recovery of
 $K$-sparse nonnegative vectors.
 This concept also enables us to develop certain conditions for the non-uniform recovery of sparse
nonnegative vectors via the so-called weak range space property.\\



\newpage

\section{Introduction}

 In this paper, we use $\|\cdot\|_0$ to denote   the number of nonzero components of a vector.
We investigate the following optimization problem with nonnegativity
constraints:
\begin{equation}\label{000} \min
\{\|x\|_0: ~ Ax=b, ~ x\geq 0\},
\end{equation} which is called an
 $\ell_0$-minimization problem or $\ell_0$-problem.
    It is well known that nonnegativity constraints are quite common in
mathematical optimization and numerical analysis (see [1] and the
references therein).  Clearly, the aim of the problem (\ref{000}) is
to find a sparsest nonnegative solution to a system of linear
equations. This problem has found so many applications in such areas
as signal and image processing [2$-$10], machine learning [11$-$15],
pattern recognition and computer vision
 [9, 16], proteomics [17],
 to name but a few. This problem
is a special case of
  compressed nonnegative sparse coding [18$-$19],
and  rank minimization  with positive semidefinite constraints
(e.g., [10, 20$-$21]). It is closely related to the  nonnegative
matrix factorization as well [22$-$24].

 The $\ell_0$-minimization problem is NP-hard [25].
Current theories and algorithms for $\ell_0$-minimization are mainly
developed through certain heuristic methods and continuous approximations. A
large amount of recent attention is attracted to the
$\ell_1$-problem
\begin{equation}\label{111} \min \{\|x\|_1: ~ Ax=b, ~ x\geq 0\}
\end{equation}  which is efficient for solving  (\ref{000}) in many
situations, so is
 the  reweighted $\ell_1$-minimization method (e.g., [26$-$27]). In this paper, the
optimal solution of the problem (\ref{111}) is called the
\emph{least $\ell_1$-norm nonnegative solution} to the linear system
$Ax=b. $ Any linear programming solver can be used to solve the
problem (\ref{111}). Various specialized algorithms for this problem
have also been   proposed  in the literature (e.g., [3, 9, 13,
28$-$29]).

Over the past few years, $\ell_0$-problems without nonnegativity
constraints have been extensively studied in the fields of sparse
signal and image processing and compressed sensing.   Both theories
and numerical methods have been developed for this problem (e.g.,
[26, 27, 30$-$34]). However, the sparsest solution and sparsest
nonnegative solution to a linear system are very different from a
mathematical point of view. The analysis and many results developed
for the sparsest solution of a linear system cannot apply to the
sparsest nonnegative solution straightaway. So far,  the
understanding of the relationship between (\ref{000}) and
(\ref{111}), and the $\ell_1$-method-based recovery theory for
sparse nonnegative vectors remains incomplete, compared with the
linear systems without nonnegativity constraints. For example, the
following important
 questions   have not
 well addressed at present:
  \begin{itemize} \item[(a)] How to completely characterize the uniqueness of least
  $\ell_1$-norm nonnegative solutions to a linear system?

  \item[(b)] Empirical results indicate that when existing sufficient criteria for the efficiency of the  $\ell_1$-method
 fail,  the  $\ell_1$-method  still succeeds  in solving $\ell_0$-problems in
many such situations. The $\ell_1$-method actually performs  remarkably
better than what the current theories have indicated. How to interpret
such a gap between the performance of the  $\ell_1$-method  indicated by the
current theories and that demonstrated by  the numerical simulations?

  \item[(c)] Are there any other matrix properties that are different
  from the existing ones (such as the restricted isometric property (RIP) [31, 35$-$36]
  and the null space property (NSP) [6, 37]),
  and  can fully characterize  the uniform recovery (i.e., the exact recovery of all $K$-sparse nonnegative vectors by a single sensing matrix)?

  \item[(d)] Is it possible to develop some theory for the recovery
  of sparse nonnegative vectors (representing signals or images) that may
go beyond the scope of the current
    uniform recovery?
\end{itemize}

In general, for a given pair $(A,b),$ the sparsest nonnegative
  solution
 to the
system $Ax=b$ is   not unique. So it is important to distinguish the
equivalence and the strong equivalence between (\ref{000}) and
(\ref{111}).  In this paper, \emph{$\ell_0$- and $\ell_1$-problems
are said to be equivalent if the $\ell_0$-problem has an optimal
solution that coincides with the unique optimal solution to the
$\ell_1$-problem. We say that the $\ell_0$- and $\ell_1$-problems
are strongly equivalent if the $\ell_0$-problem has a unique optimal
solution that coincides with the unique optimal solution to the
$\ell_1$-problem.} Clearly, the `strong equivalence' implies the
`equivalence', but the converse is not true in general.
 The `equivalence' does not require an $\ell_0$-problem to have a
unique optimal solution.  As shown by our later analysis,
the `equivalence' concept enables us to   broadly
understand the relationship between $\ell_0$- and
$\ell_1$-minimization, making it possible to address the
aforementioned questions.
  Of course, the above-mentioned questions (a)-(d) can be
partially addressed by applying the existing theories
 based on such concepts as the mutual coherence  [38$-$40, 45], ERC
[41$-$42], RIP [31, 35$-$36], NSP [6, 37], outwardly
$k$-neighborliness property [4], and the verifiable condition
[43$-$44]. However,   these existing conditions are relatively restrictive in
the sense that they imply the strong equivalence (instead of the
equivalence) between $\ell_0$- and $\ell_1$-problems. For instance,
Donoho and Tanner [4] have given a  geometric condition, i.e., the
outwardly $K$-neighborliness property of a sensing matrix, which
guarantees that if a $K$-sparse nonnegative vector is a solution to the linear system $Ax=b,$ then it must be the  unique optimal solution
 to  both  problems (\ref{000}) and (\ref{111}). From a
null-space perspective, Khajehnejad et al [6] have shown that
$K$-sparse nonnegative vectors can be recovered by
$\ell_1$-minimization if and only if the null space of $A$ satisfies
certain property. Thus both the outwardly $K$-neighborliness property [4]
and the null space property [6] imply the strong equivalence between
 (\ref{000}) and (\ref{111}). We note that  the strong equivalence conditions fail to explain the success of the $\ell_1$-method  for solving
 $\ell_0$-problems with  multiple optimal solutions. We also note that the uniqueness of least $\ell_1$-norm nonnegative solutions to a linear system   plays a fundamental role in both theoretical and practical
 efficiencies of the $\ell_1$-method  for solving  $\ell_0$-problems.  While the  strong equivalence conditions   are
sufficient for the uniqueness of least $\ell_1$-norm nonnegative solutions to a linear system,  these
conditions are not the necessary condition.

The first purpose of this paper is to completely address the
question (a) by developing  a necessary and sufficient condition
 for the uniqueness of least
$\ell_1$-norm nonnegative solutions to a linear system. We establish
this  condition through the   strict complementarity theory of
linear programming, which  leads naturally to the new concept of the
range space property (RSP) of $A^T. $    Based on this result, we
further point out that the equivalence between $\ell_0$- and
$\ell_1$-problems can be theoretically interpreted by   the RSP of $A^T
.$  That is,  an $\ell_1$-problem is equivalent  to   an
$\ell_0$-problem if and only if the RSP of $A^T$ holds at an optimal solution
of the $\ell_0$-problem.

While the RSP of $A^T$ is  defined locally  at an individual vector (e.g., the solution to an $\ell_1$- or an $\ell_0$-problem),
 it provides a complete and practical checking condition for the uniqueness of least  $\ell_1$-norm nonnegative solutions to a linear system, and
the RSP of $A^T$ yields a
broad understanding of the efficiency  of  the
$\ell_1$-method for solving $\ell_0$-problems (see the discussion in Sect. 3 for details). It turns out that when the
strong equivalence conditions fail (e.g., the $\ell_0$-problem has multiple optimal solutions),  the RSP of $A^T$ is still
 able to explain the success of the $\ell_1$-method  for  solving  $\ell_0$-problems in many such  situations. Thus,
 the current gap between the   performance of the $\ell_1$-method indicated by  the  existing theories and that observed from numerical
simulations can be clarified  in terms  of the RSP of $A^T, $ leading to a certain answer to the question (b).

Although a global RSP-type condition for the equivalence between $\ell_0$- and $\ell_1$-problems
remains not clear at present, such a condition  for the strong
equivalence between  these two problems can be developed.     In Sect. 4,  we further introduce
 a matrix property, called the RSP of
order $K$, through which we provide a characterization of the uniform
recovery of all $K$-sparse nonnegative vectors.   As a result, the RSP of order $K$ is a
strong equivalence condition for $\ell_0$- and $\ell_1$-problems.   Interestingly, the
variants of this new concept make it possible to extend uniform
recovery to  non-uniform recovery of some sparse nonnegative
vectors, to which the uniform recovery does not apply. Such an
extension is important not only from a mathematical point of view,
but from the viewpoint of many practical applications as well. For
instance, when many columns of $A$ are important, the sparsest
solution to the linear system $Ax=b$ may not be sparse enough to
satisfy  the uniform recovery condition. The RSP of order $K$ and
its variants make it possible to address the aforementioned
questions (c) and (d).

This paper is organized as follows. In Sect. 2, we develop a
necessary and sufficient condition for a linear system to have a
unique least $\ell_1$-norm nonnegative solution.
  In Sect. 3, we provide an efficiency analysis for the  $\ell_1$-method
   through the RSP of $A^T.$ In Sect. 4, we
  develop some (uniform and non-uniform)
recovery conditions for $K$-sparse nonnegative vectors via the so-called RSP of
order $K$ and its variants.  Conclusions are given in the last section.

\section{Uniqueness of least $\ell_1$-norm nonnegative solutions}

We use the following notation:   Let $ R^n_{+}$
be the first orthant of $R^n,$ the $n$-dimensional Euclidean space.
Let $e= (1,1,...,1)^T\in R^n$ be the vector of ones throughout this
paper. For two vectors $u,v\in R^n$, $u\leq v$ means $u_i \leq v_i $
for every $i=1, ..., n, $ and in particular, $v\geq 0$ means $ v\in
R_+^n.$ For a set $S\subseteq \{1,2,..., n\},$ $|S|$ denotes the
cardinality of $S,$ and $S_c = \{1, 2,..., n\} \backslash S $ is the
complement of $S.$ For a matrix $A$ with columns $a_j,$ $1\leq j\leq
n,$ we use $A_S$ to denote the submatrix of $A$ with columns $a_j,
j\in S.$ Similarly, $x_S $ denotes the subvector of $x$ with
components $x_j,$  $j\in S.$ For $x\in R^n,$ let
$\|x\|_1=\sum_{i=1}^n |x_j| $ denote the $\ell_1 $-norm of $x.$ For
$A\in R^{m\times n}, $ we use ${\cal R} (A^T) $ to denote the range
space of $A^T, $ i.e., $ {\cal R} (A^T) = \{A^T u: u \in R^m\}.$

In this section, we develop a necessary and sufficient condition for
$x$ to be the unique least $\ell_1$-norm nonnegative solution to a
linear system.
 Note that when $x$ is the unique
 optimal solution to the problem (\ref{111}),  there is no other nonnegative solution
 $w\not= x$ such that  $\|w\|_1\leq \|x\|_1.$   Thus the uniqueness of  $x$ is equivalent to
$$ \{w: ~ Aw=b,~  w\geq 0, ~\|w\|_1 \leq \|x\|_1\}  =\{x\}.
$$
Since $x\geq 0$ and $ w\geq 0,$ we have $\|w\|_1=e^Tw $ and $\|x\|_1
=e^T x. $  Thus the above relation can be further written as  $ \{w:
~ Aw=Ax,~ e^T w \leq e^T x,  ~  w \geq 0 \}  =\{x\}.  $ Consider the
following linear programming (LP) problem  with the variable $w\in
R^n:$
\begin{equation} \label{LP1}  \min \{0^T w: ~Aw=Ax,  ~ e^T w \leq
e^T x, ~ w \geq 0 \},
\end{equation}  which is feasible (since $w=x$ is always a feasible solution),
and the optimal value of the problem is finite (equal to zero). From
the above discussion, we immediately have the following observation.

\vskip 0.08in

\textbf{Lemma 2.1 } \emph{ $x$ is  the unique least $\ell_1$-norm
nonnegative solution to the system $Ax=b $  if and only if $w=x$ is
the unique optimal solution to the problem (\ref{LP1}), i.e., $(w,t)
= (x,0)$ is the unique optimal solution to the following problem:
\begin{equation} \label{LP2} \min\{ 0^T w: ~Aw=Ax,  ~ e^T w +t= e^T x, ~ (w,t) \geq 0
\}\end{equation}
 where $t$ is a slack variable introduced into (\ref{LP1}). }

\vskip 0.08in

Note that the dual problem of (\ref{LP2}) is given by
\begin{equation}  \label{DLP2}     \max   \{(Ax)^T y+(e^Tx)\beta: ~
                     A^T y +\beta e \leq 0,   ~ \beta \leq 0  \} \end{equation}
  where $y$ and $\beta$ are variables.
  Throughout this section, we use $(s, r)\in
  R^{n+1}_+$ to denote the slack variables of the
  problem (\ref{DLP2}), i.e.,
$$   s= - (A^T y +\beta e) \geq 0, ~~ r= -\beta \geq
 0.$$

Let us   recall a fundamental theorem for LP problems. Let $B\in
R^{m\times n}$ be a given  matrix, and  $p\in R^m$ and $  c\in R^n$
be two given vectors.  Consider the  LP problem
\begin{equation} \label {LP}    \min \{c^T x: ~ Bx=  p, ~~ x\geq 0\},
\end{equation}
and its dual problem
\begin{equation} \label{DP}
\max\{p^Ty: ~ B^T y + s= c, ~ s\geq 0\}. \end{equation} Any optimal
solution pair $(x, (y,s))$ to  (\ref{LP}) and (\ref{DP})
satisfies the so-called complementary slackness condition:  $x^Ts=0,
x\geq 0 $ and $ s\geq 0. $ Moreover, if a solution pair $(x, (y,
s))$ satisfies that $x+s>0,$  it is called  a strictly complementary
solution pair. For any feasible LP problems
(\ref{LP}) and (\ref{DP}), there always exists a pair of strictly
complementary solutions.

\vskip 0.08in

 \textbf{Lemma 2.2}([46]) ~  (i) (Optimality condition)  \emph{$(x, (y, s))$ is a solution pair of the
 LP problems (\ref{LP}) and (\ref{DP}) if and only if it
satisfies the following conditions: $
                             Bx =  p,  ~B^T y +s = c , ~~ x\geq 0, ~ s\geq 0, $ and $
                             x^Ts    =   0.$ }   (ii) (Strict complementarity)
                             \emph{ If (\ref{LP}) and (\ref{DP}) are feasible,  then there exists a pair
$(x^*, (y^*, s^*))$ of strictly complementary solutions to
(\ref{LP}) and (\ref{DP}). }

\vskip 0.08in

We now prove the following   necessary condition for the problem
(\ref{111}) to have a  unique optimal solution.

\vskip 0.08in

\textbf{Lemma 2.3 }  \emph{If $x$ is the unique least $\ell_1$-norm
nonnegative solution to the system $Ax=b, $ then there exists a
vector $\eta\in R^n$ satisfying \begin{equation} \label{RSP+} \eta
\in { \mathcal R} (A^T), ~ \eta_i =1  \textrm{ for } i\in J_+,
\textrm{ and } \eta_i <1  \textrm{ for } i\notin J_+,
\end{equation}
 where $J_+ = \{ i: x_i >0\}. $  }

\vskip 0.08in

  \emph{Proof.}  Consider the problem (\ref{LP2}) and its dual problem (\ref{DLP2}), both of
 which are feasible.   By Lemma 2.2,    there exists an optimal solution $ (w^*, t^*)$
 to the problem
(\ref{LP2})  and  an optimal solution $(y^*, \beta^*)$ to
(\ref{DLP2}) such that these two solutions constitute a pair of
strictly complementary solutions. Let $(s^*, r^*) =  ( - A^T y^*
-\beta^* e, -\beta^*)
 $ be the value of the associated  slack  variables of the dual problem (\ref{DLP2}). Then by the strict complementarity, we have
\begin{equation} \label {strict-comp}  ( w^*)^T s^* =0, ~ t^* r^* =0,  ~ w^*+s^* >0, ~ t^*+ r^*
>0.
\end{equation}
Since $x $ is the unique least $\ell_1$-norm nonnegative solution to
$Ax=b,$ by Lemma 2.1, $(x, 0)$ is the unique optimal solution to the
problem (\ref{LP2}). Thus
\begin{equation} \label{www} (w^*, t^* ) = (x,  0), \end{equation}
which  implies that $w^*_i > 0$ for all $i\in J_+ =: \{i: x_i>0\}$
and $w^*_i=0$ for all $i\notin J_{+}.$ Thus it follows from
(\ref{strict-comp}) and (\ref{www}) that $r^*>0, ~  s^*_i =0~
\textrm{ for all } i\in J_+, $   and that
$ s^*_i >0 $ for all $ i\notin J_+.  $
That is,
$$\beta^* <0, ~ (A^Ty^*+\beta^* e)_i =0 ~ \textrm{ for  } i\in J_+,
~ (A^Ty^*+\beta^* e)_i <0 ~ \textrm{ for  } i\notin J_+, $$ which
can be written as
$$\beta^* <0,  ~ \left[A^T(\frac{y^*}{-\beta^*})- e\right]_i =0 ~ \textrm{ for  } i\in J_+,
~ \left[A^T( \frac{y^*}{-\beta^*})-e\right]_i <0 ~ \textrm{ for  }
i\notin J_+.
$$ By setting $\eta= A^Ty^*/(-\beta^*), $ the condition above
is equivalent to
$$ \eta \in {\mathcal R }(A^T), ~ \eta_i=1  ~ \textrm{ for  } i\in J_+, \textrm{  and } \eta_i
<1 ~ \textrm{ for } i\notin J_+, $$ as desired. \hfill  $\Box$

Throughout this paper, the condition (\ref{RSP+}) is called the
\emph{range space property }  (RSP) of $A^T$ at $x\geq 0.$  It is worth noting  that   the RSP (\ref{RSP+}) can be easily checked   by simply solving the following LP problem:
  $$ t^*= \min\left\{t:  ~(A_{J_+})^T y= e_{J_+},  ~(A_{J_{+c}})^Ty = \eta_{J_{+c}}, ~\eta_{J_{+c}} \leq t e_{J_{+c}} \right\}, $$  where $ J_{+c}=\{1,2, ..., n\}\setminus  J_+. $ Clearly,  $t^*<1$ if and only if the RSP (\ref{RSP+}) holds.  Lemma
2.3 shows that the RSP of $A^T$ at $x$
 is a necessary condition for the $\ell_1$-problem to have a unique optimal solution.  We now
prove  another necessary condition.

\vskip 0.08in

 \textbf{Lemma 2.4 }  \emph{ If
$x$ is the unique least $\ell_1$-norm nonnegative solution to the
system $Ax=b, $ then the matrix \begin{equation} \label{MMM} M=
\left(
                              \begin{array}{c}
                                A_{J_+} \\
                                e^T_{J_+} \\
                              \end{array}
                            \right) \end{equation}  has full column rank,
 where $J_+ = \{ i: x_i >0\}. $ }

\vskip 0.08in

 \emph{Proof.} Assume the contrary that the columns of  $M $ defined by (\ref{MMM})
 is linearly dependent. Then there exists a vector
                            $u\in R^{|J_+|} $ such that
                        \begin{equation}\label{uuu} u\not =0, ~   Mu=\left(
                              \begin{array}{c}
                                A_{J_+} \\
                                e^T_{J_+} \\
                              \end{array}
                            \right)    u=0.
                            \end{equation}
Let $(w,t)$ be given by  $w= (w_{J_+}, w_{J_0}) = (x_{J_+}, 0)$ and
$t=0, $ where $J_0=\{i: i\not\in J_+\}.$  Then it is easy to see
that such defined $(w,t)$ is an optimal solution to the problem
(\ref{LP2}). On the other hand, let us define
$(\widetilde{w},\widetilde{t})$ as follows:
$$ \widetilde{w} = (\widetilde{w}_{J_+}, \widetilde{w}_{J_0}) =(
w_{J_+}+\lambda u, 0),  ~  \widetilde{t} =0.$$ Since $w_{J_+}
=x_{J_+}
>0, $ there exists a small $\lambda \not=0 $ such that
\begin{equation} \label{YYY} \widetilde{w}_{J_+} = w_{J_+}+\lambda u \geq 0.\end{equation} Substituting
$ (\widetilde{w},\widetilde{t})$ into the constraints of the problem
(\ref{LP2}),
  we see from (\ref{uuu}) that $(\widetilde{w},\widetilde{t})$
satisfies all those constraints. Thus
$(\widetilde{w},\widetilde{t})$ is also an optimal solution to the
problem (\ref{LP2}). It follows from (\ref{YYY}) that $
\widetilde{w}_{J_+} \not= w_{J_+}$ since $\lambda u\not=0.$
Therefore, the optimal solution to   (\ref{LP2}) is not
unique. However, by Lemma 2.1, when $x$ is the unique least
$\ell_1$-norm nonnegative solution to the system $Ax=b,$ the problem
(\ref{LP2}) must have  a unique optimal solution. This contradiction
shows that  $M $ must have full column rank. \hfill $\Box$

    The next result shows that the combination of the necessary
conditions developed  in Lemmas 2.3 and 2.4 is  sufficient  for  the
$\ell_1$-problem to have a unique optimal solution.

\vskip 0.1in

\textbf{Lemma 2.5 } \emph{Let $x\geq 0 $ be a solution to the system
$Ax=b.$
  If the condition (\ref{RSP+})  (i.e., the RSP of $A^T$) is
satisfied at $x,$ and  the  matrix
      $M$ given by (\ref{MMM}) has full column rank,  then $x$ is the unique least $\ell_1$-norm
nonnegative solution to the system  $Ax=b.$}

\vskip 0.08in

\emph{Proof.} By Lemma 2.1, to prove that $x $ is the unique least
$\ell_1$-norm nonnegative solution to the system $Ax=b,$ it is
sufficient to prove that the problem (\ref{LP2}) has a unique
optimal solution $(x,0).$ First, the condition (\ref{RSP+}) implies
that there exist $ \eta $ and $y$ such that $$ A^T y=\eta, ~ \eta_i
=1 \textrm{  for } i\in J_+, \textrm{ and } \eta_i<1 \textrm{ for
}i\notin J_+.$$ By setting $\beta =-1$, the  relation above can be
written as
\begin{equation} \label{dual-sol} (A^T y)_i +\beta =0 ~ \textrm{ for }i\in J_+, ~ \textrm{ and }
(A^T y)_i +\beta < 0   \textrm{ for } i\notin J_+,\end{equation}
from which  we see that $(y, \beta)$   satisfies all constraints of
the  problem (\ref{DLP2}).   We now further  verify that it is an
optimal solution to   (\ref{DLP2}). In fact, by (\ref{dual-sol}), the
objective value of (\ref{DLP2}) at $(y, \beta)$ is
\begin{eqnarray} (Ax)^T y +
(e^Tx) \beta  & = &  x^T (A^T y) + (e^Tx) \beta  \nonumber\\
& =  & \sum_{i\in J_+} x_i (A^Ty)_i + \beta \sum_{i\in J_+}x_i \nonumber\\
 & =  & -\beta\sum_{i\in J_+} x_i + \beta \sum_{i\in J_+}x_i =0.
\label{444}
\end{eqnarray}
Since the optimal value of (\ref{LP2}) is zero, by LP duality
theory, the maximum value of the dual problem is also zero. Thus it
follows from (\ref{444}) that the point $(y, \beta)$ satisfying
(\ref{dual-sol}) is an optimal solution to the  problem
(\ref{DLP2}).

We now prove that the optimal solution of (\ref{LP2}) is uniquely
determined under the assumption of the theorem. Assume that
$(w^*,t^*)$ is an arbitrary optimal solution to the problem
(\ref{LP2}),
 which of course satisfies all constraints of (\ref{LP2}), i.e.,
\begin{equation} \label{LP-constr}  Aw^*=Ax,  ~ e^T w^* +t^*= e^T x,
~ (w^*,t^*) \geq 0. \end{equation}   Since  $(y, \beta)$ satisfying
(\ref{dual-sol}) is an optimal solution of (\ref{DLP2}), $((w^*,
t^*), ((y, \beta),s))$ is a solution pair to (\ref{LP2}) and
(\ref{DLP2}). From (\ref{dual-sol}), we see that the dual slack
variables $ s_i = -((A^T y)_i+\beta)
>0$ for $ i\notin J_+ $ and $r =-\beta =1 >0. $  By complementary slackness property (Lemma 2.2(i)), we
must have that
$$ t^*=0,  ~ w^*_i =0\textrm{ for all }  i\notin J_+.  $$
By substituting these known components into (\ref{LP-constr}) and
noting that $x_{i}=0 $ for $i\notin J_+,$   we see that the
remaining components of $(w^*, t^*)$ satisfy
$$A_{J_+} w^*_{J_+} = Ax =A_{J_+} x_{J_+},  ~ e^T_{J_+}w^*_{J_+} = e^T x =e^T_{J_+} x_{J_+},
~ w^*_{J_+} \geq 0 . $$ Since
 the matrix {\small $M=\left(
     \begin{array}{c}
       A_{J_+}      \\
       e^T_{J_+}
     \end{array}
   \right)
$}  has full column rank,   $w^*_{J_+}= x_{J_+} $ is the unique
solution to the  reduced system above. Therefore, $(w^*, t^*) $ is
uniquely given by $(x,0). $ In other words, the only optimal
solution to the problem (\ref{LP2}) is $(x, 0). $
 By Lemma 2.1, $x$ must be the unique
least $\ell_1$-norm  nonnegative solution to the system $Ax=b.$
\hfill $\Box$

    By Lemmas 2.3, 2.4 and 2.5, we have the main result of this
section.

\vskip 0.08in

 \textbf{Theorem 2.6 } \emph{Let $x$ be a nonnegative solution to the system $Ax=b.$ Then $x$
is the unique  least $\ell_1$-norm nonnegative solution to the
system $Ax=b $ if and only if the RSP (\ref{RSP+}) holds at $x$ and
the matrix $M=\left(
     \begin{array}{c}
       A_{J_+}      \\
       e^T_{J_+}
     \end{array}
   \right)
$ has full column rank, where $J_+=\{i: x_i>0\}. $ }

\vskip 0.08in

 Clearly, when $A_{J_+}$ has full column
   rank, so does the matrix $M$ given by (\ref{MMM}).
The converse is not true.  In general, when   $M $
 has full column rank, it does not imply that the matrix $A_{J_+}$
has full column rank. For instance, {\small M=  $
\left(\begin{array}{c}
       A_{J_+}      \\
       e^T_{J_+}
     \end{array}
   \right)= \left(
              \begin{array}{cc}
                -1& 1 \\
                0 & 0 \\
                1 & 1 \\
              \end{array}
  \right)$ }  has full column rank, but {\small $A_{J_+} = \left(
\begin{array}{cc}
 -1 & 1 \\
  0 & 0 \\
  \end{array}
   \right)$ } does not. However, \emph{when the RSP (\ref{RSP+}) holds at $x,$
   we see that $e_{J_+} = A^T_{J_+} u $ for some $u\in R^n, $ in
   which case $A_{J_+}$ has full column rank if and only if $\left(\begin{array}{c}
       A_{J_+}      \\
       e^T_{J_+}
     \end{array}
   \right)$ has full column rank.} Thus Theorem 2.6 can be restated
   as follows.

   \vskip 0.08in

 \textbf{Theorem 2.7 } \emph{Let $x$ be a nonnegative solution to the system $Ax=b.$ Then $x$
is the unique  least $\ell_1$-norm nonnegative solution to the
system $Ax=b $ if and only if the RSP (\ref{RSP+}) holds at $x$ and
the matrix $ A_{J_+} $ has full column rank, where $J_+=\{i:
x_i>0\}. $   }

\vskip 0.08in

The above results   completely characterize the uniqueness of least
$\ell_1$-norm nonnegative solutions to a system of linear equations,
and thus the question (a) in Sect. 1 has been fully addressed. Note
that $A\in R^{m\times |J_+|}, $ so when it has full column rank, we
must
   have $\textrm{rank}(A_{J_+}) =|J_+|
   \leq m.$  Thus Theorem 2.7 shows that \emph{if the $\ell_1$-problem  has a unique optimal solution $x$,
   then $x$  must be  $m$-sparse.}
We can use the results established in this section to discuss
other questions associated with $\ell_0$- and $\ell_1$-problems (see the remainder of the paper for details).

We now close this section by giving two
 examples to show that our necessary and sufficient condition
 can be easily used to check the uniqueness of least
$\ell_1$-norm nonnegative solutions of linear systems.

\vskip 0.08in
  \textbf{Example 2.8 } Consider the linear
system $Ax=b$  with $$A= \left(
       \begin{array}{cccc}
         1 & 0 & -1 & -1 \\
         0 & -1 & -1 & 6 \\
         0 & 0 & -1 & 1 \\
       \end{array}
     \right)
 , ~~  b=  \left(
          \begin{array}{c}
            1/2 \\
            -1/2 \\
            0 \\
          \end{array}
        \right),
  $$ to which   $x^*= (1/2,1/2, 0, 0)^T$
  is a nonnegative solution. It is easy to see that  the submatrix $A_{J_+}$
   associated with this solution has full column rank.
     Moreover, by taking $y= (1, -1,
    0)^T$, we have $\eta= A^T y =  (1, 1, 0, -7)^T  \in {\cal R}(A^T), $  which clearly satisfies (\ref{RSP+}).
     Thus
    the RSP of $A^T$  holds at $x^*.$ Therefore, by Theorem  2.7 (or Lemma 2.5), $x^*$ is the unique least $\ell_1$-norm nonnegative   solution
     to the  system
    $Ax=b.$

 \vskip 0.08in  \textbf{Example 2.9 } Consider the linear
system $Ax=b $ with $$A= \left(
       \begin{array}{crrr}
         1 & 0 & -1 & 1 \\
         1 & -0.1 & 0 & -0.2 \\
         0 & 0 & -1 & 1 \\
       \end{array}
     \right)
 , ~~  b=  \left(
          \begin{array}{r}
            1/2 \\
            -1/2 \\
            0 \\
          \end{array}
        \right),
  $$ to which  $x^*= (1/2,10/3, 10/3, 10/3)^T$
    is a least $\ell_1$-norm nonnegative solution.  By taking $y= (11, -10,
    -12)^T$, we have $\eta= A^T y =  (1, 1, 1,1)^T  \in {\cal R}(A^T).$ Thus
    the RSP of $A^T$ holds at $x^*.$  However, the matrix  $$
                           A_{J_+}
     = \left(\begin{array}{crrr}
         1 & 0 & -1 & 1 \\
         1 & -0.1 & 0 & -0.2 \\
         0 & 0 & -1 & 1
       \end{array}
     \right)
     $$  does not have full column rank.
      By Theorem  2.7, $x^*$ is not the unique least $\ell_1$-norm nonnegative solution
     to the  system
    $Ax=b.$ In fact, we have another least $\ell_1$-norm nonnegative solution given by $\widetilde{x} =(1/2,10, 0,
    0)^T$ (for which  the associated $ A_{J_+}$
    has full column rank, but the RSP of $A^T$ does not hold at
    $\widetilde{x}$).

\section{RSP-based efficiency analysis for $\ell_1$-minimization}
  For
linear systems without nonnegativity constraints, some sufficient
conditions  for the strong equivalence between $\ell_0$- and
$\ell_1$-problems have been developed in the literature.
 If  these sufficient conditions are applied directly to  sparsest
 nonnegative  solutions of linear systems, the resulting criteria would be very restrictive.
  For instance, by applying the mutual
coherence condition, we immediately conclude that if a nonnegative
solution $x$ obeys $\|x\|_0<(1+1/\mu(A))/2$ where $\mu(A)$ denotes
the mutual coherence of $A$ (i.e., $\mu(A) =\max_{i\not= j} |a_i^T
a_j|/(\|a_i\|_2\|a_j\|_2) $ where $a_i, 1\leq i \leq n,$ are the
columns of $A$),  then $x$ is the unique sparsest solution and the
unique least $\ell_1$-norm solution to the linear system $Ax=b.$ In
this case, the unique sparsest nonnegative solution coincides with
the unique sparsest solution and the unique least $\ell_1$-norm
solution of the linear system. Clearly, such a sufficient  (strong-equivalence-type) condition
is too restrictive. In fact, a sparsest nonnegative solution is
usually not a sparsest one of the linear system,  and the sparsest
nonnegative ones can be   multiple (as shown by Example 3.4 in
this section). Although some conditions have been developed
specifically for the sparsest nonnegative solution in the literature (e.g., [3$-$4,
6]),
 these conditions  still   imply the strong equivalence between $\ell_0$- and
 $\ell_1$-problems.  They can only partially explain the success of the
 $\ell_1$-method for solving $\ell_0$-problems.   In this
 section, we point out that Theorem 2.6 or 2.7  can be used  to broadly clarify the relationship
 between
 $\ell_0$- and $\ell_1$-problems,
 and thus the efficiency of the $\ell_1$-method can be more widely interpreted
 through the
RSP of $A^T.$ First,  we have the following property for  sparsest
nonnegative solutions.

\vskip 0.08in

 \textbf{Lemma 3.1 } \emph{ If $x$ is a sparsest nonnegative solution to the system  $Ax=b$, then {\small $M = \left(
     \begin{array}{c}
       A_{J_{+}}      \\
       e^T_{J_{+}}
     \end{array}
   \right)
$ } has full column rank, where $J_+=\{i: x_i>0\}.$   }

\vskip 0.1in

\emph{Proof.}  Let $x$ be a sparsest nonnegative solution to the
system $Ax=b$ and let $J_+=\{i: x_i>0\}.$ Assume by contrary that
the columns of   $M$ are linearly dependent. Then there
exists a vector $v \not =0$ in $R^{|J_+|}$ such that { \small $$
\left(
     \begin{array}{c}
       A_{J_{+}}      \\
       e^T_{J_{+}}
     \end{array}
   \right) v=0.$$ }
It follows from  $e^T_{J_{+}}v=0$ and $v\not=0$  that $v$ must have
at least two nonzero components with different signs, i.e., $v_iv_j
<0 $ for some $i\not= j. $ Define the vector $ \widetilde{v} \in R^n
$ as follows: $\widetilde{v}_{J_+} = v $  and $\widetilde{v}_i=0$
for all $i\notin J_+ .$ We consider the vector $y(\lambda)=
x+\lambda \widetilde{v}$ where $ \lambda \geq 0.$ Note that
$y(\lambda)_i =0$ for all $i\notin J_+, $ and that
$$ Ay(\lambda) = Ax+A (\lambda \widetilde{v}) =b +  \lambda A_{J_+} v =b.$$
Thus $y(\lambda)$ is also a solution to the linear system $Ax=b. $
By the definition of $\widetilde{v},$ $ \widetilde{v}$ has  at least
one negative component. Thus let
$$\lambda^* = \frac{x_{i_0}}{-\widetilde{v}_{i_0}} =\min
\left\{\frac{x_i}{-\widetilde{v}_i}: ~ \widetilde{v}_i <0\right\},
$$ where $\lambda^*$ must be a positive number and $i_0\in J_+.$  By
such a choice of $\lambda^*$ and the definition of $y(\lambda^*),$
we conclude that $y(\lambda^*) \geq 0$,  $y(\lambda^*)_i=0$ for
$i\not\in J_+,$ and $y(\lambda^*)_{i_0}=0$ with $i_0\in J_+.$ Thus
$y(\lambda^*)$ is a nonnegative solution
 to the linear system $Ax=b,$  which is sparser than $x.$ This is a contradiction.
 Therefore,  $M$ must have full column
rank. \hfill $\Box$

 By Theorem 2.6 and Lemma 3.1, we  immediately have
the following  result.

\vskip 0.08in

 \textbf{Theorem 3.2 }  \emph{$\ell_0$- and
$\ell_1$-problems are equivalent if and only if the RSP (\ref{RSP+})
holds at an optimal solution of the $\ell_0$-problem. (In other
words, a sparsest nonnegative solution $x$ to the system $Ax=b$ is
the unique least $\ell_1$-norm nonnegative solution to the system if
and only if the RSP (\ref{RSP+}) holds at $x.$) }

\vskip 0.08in

 \emph{Proof.} Assume that  problems (\ref{000}) and (\ref{111}) are equivalent.
  So the $\ell_0$-problem has an optimal solution $x$ that is the unique least  $\ell_1$-norm nonnegative
  solution to the
  system $Ax=b.$ Thus, by Theorem 2.6 (or Lemma 2.3),  the RSP
(\ref{RSP+}) must hold at $x.$  Conversely,  assume that the RSP
(\ref{RSP+}) holds at an optimal solution  $x$ to the
$\ell_0$-problem. Since $x$ is a sparsest nonnegative solution to
the system $Ax=b,$ by Lemma 3.1, the matrix {\small $\left(
     \begin{array}{c}
       A_{J_{+}}      \\
       e^T_{J_{+}}
     \end{array}
   \right)$ }   has full column rank. Thus  by Lemma 2.5 (or Theorem 2.6) again,
    $x$ must be
    the unique least
   $\ell_1$-norm nonnegative solution to the system $Ax=b.$ Hence $\ell_0$- and $\ell_1$-problems
    are equivalent. \hfill $\Box$

 It should be pointed out that the
equivalence   between $\ell_0$- and $\ell_1$-problems are characterized implicitly  in Theorem 3.2    in the sense that the RSP condition therein is defined locally
at a  solution of the $\ell_0$-problem.
 Whether or not a checkable   RSP-type equivalence (instead of strong equivalence) condition can be developed  for $\ell_0$- and $ \ell_1$-problems
  remains an open question.  However,   Theorem 3.2 is still important from a theoretical point of view, and it can be used to explain
 the numerical performance  of the $\ell_1$-method more efficiently
than strong equivalence conditions.
Since the RSP (\ref{RSP+}) at an optimal solution of $\ell_0$-problem is a necessary and
sufficient condition for the equivalence between $\ell_0$- and
$\ell_1$-problems,  all existing sufficient conditions for
 strong  equivalence (or  equivalence) between these two problems must
imply the RSP (\ref{RSP+}). However,  the converse is clearly not true in
general, as shown by the following example.

\vskip 0.08in

 \textbf{Example 3.3} (When  existing criteria fail, the RSP may still
 succeed).
 $$A= \left(
   \begin{array}{cccccc}
     0 & -1 & \frac{1}{\sqrt{3}} & 0 & \frac{1}{\sqrt{2}} & -\frac{1}{\sqrt{2}} \\
     0 & 0 & \frac{1}{\sqrt{3}} & -1 & 0 & 0 \\
     -1 & 0 & \frac{1}{\sqrt{3}} & 0 &  \frac{1}{\sqrt{2}}& -\frac{1}{\sqrt{2}} \\
   \end{array}
 \right), ~ b= \left(
               \begin{array}{c}
                 1 \\
                 1 \\
                 0 \\
               \end{array}
             \right).
$$
 For this example, the system $Ax=b$ does not have a
solution $x$ with $\|x\|_0=1. $ So $ x^*= (1, 0, \sqrt{3}, 0, 0,
0)^T$ is a sparsest nonnegative solution of this linear system.
  Note that the mutual
coherence $\mu(A)= \max_{i\not=j} |a_i^Ta_j|/\|a_i\|_2 \|a_j\|_2=
\sqrt{2}/\sqrt{3}.$ Thus the mutual coherence condition $\|x\|_0 <
\frac{1}{2}(1+1/\mu(A))= (\sqrt{2}+\sqrt{3})/(2\sqrt{2}) \approx
1.077$ fails for this example.  The RIP [35]  fails since the last
two columns of $A$ are linearly dependent. This example also fails
to comply with the definition of the NSP.  Let us now check the RSP
of $A^T$ at $x^*.$ By taking $y= (\frac{1}{2}+\sqrt{3}, \frac{1}{2},
-1)^T$, we have
$$\eta = A^T y = \left(1, -(\frac{1}{2}+\sqrt{3}), 1, -\frac{1}{2},
\frac{2\sqrt{3}-1}{2\sqrt{2}}, -
\frac{2\sqrt{3}-1}{2\sqrt{2}}\right)^T \in {\cal R }(A^T),$$ where
the first and third components of $\eta$ are equal to 1
(corresponding to $J_+=\{1,3\}$ determined by $x^*$)  and all other
components of $\eta$ are less than 1. Thus the RSP (\ref{RSP+})
holds at $x^*.$ By Theorem 3.2,  $\ell_1$-minimization will find
  this solution.

This example indicates that even if the existing sufficient
 conditions fail,  the RSP of $A^T$ at a vector  may
 still be able to confirm the success of  the $\ell_1$-method  when solving
an $\ell_0$-problem.
To further understand the efficiency of  the $\ell_1$-method, let us
decompose the class of linear systems  with nonnegative solutions,
denoted by ${\cal G} $, into three subclasses. That is, ${\cal G} =
{\cal G}_1 \bigcup {\cal G}_2 \bigcup {\cal G}_3$ where ${\cal
G}_i$'s
 are defined as follows:

\begin{enumerate}
  \item [${\mathcal G}_1$]   The  system $Ax=b$  has a \emph{unique} least $\ell_1$-norm nonnegative solution
  and a \emph{unique} sparsest nonnegative solution.

 \item [${\mathcal G}_2$]    The  system $Ax=b$  has a \emph{unique} least $\ell_1$-norm nonnegative
solution and \emph{multiple} sparsest  nonnegative  solutions.

  \item [${\mathcal G}_3$]  The  system $Ax=b$  has \emph{multiple} least $\ell_1$-norm
  nonnegative solutions.
\end{enumerate}
Clearly, every linear system   with a nonnegative solution falls
into one of these categories. Since many existing sufficient
conditions (such as the mutual coherence, RIP and NSP) imply the
strong equivalence between $\ell_0$- and $\ell_1$-problems, these
  conditions can apply only to (and explain the efficiency of the $\ell_1$-method only for) a subclass of linear systems in
${\mathcal G}_1. $  However, the RSP (\ref{RSP+}) defined in this
paper goes   beyond this scope of linear systems.  An important
feature of the RSP (\ref{RSP+}) is that it does not require a linear
system to have a unique sparsest nonnegative solution  in order to
achieve the equivalence between $\ell_0$- and $\ell_1$-problems, as
shown by the next example.

\vskip 0.08in

 \textbf{Example 3.4 } (The $\ell_1$-method may successfully solve $\ell_0$-problems with multiple
 optimal solutions.)
  Consider the system $Ax=b$ with
 $$ A=  \left(
       \begin{array}{rrrrrr}
         0.2 & 0 & -0.3 & -0.1 & 0.5 & -0.25 \\
         0 & 0.2 & 0.5 & 0.2 & -0.9 & 0.05 \\
         0.2 & 0 & -0.3 & -0.1 & 0.5 & -0.25 \\
       \end{array}
     \right), ~ b= \left(
                     \begin{array}{r}
                       0.1 \\
                       -0.1 \\
                       0.1 \\
                     \end{array}
                   \right).
 $$
For this example, it is easy to verify that $Ax=b$ has multiple
sparsest nonnegative solutions:
$$x^{(1)} = (0, \frac{2}{5}, 0, 0, \frac{1}{5}, 0)^T,~  x^{(2)} = (0, 0, 0, 4,  1, 0)^T,
~x^{(3)} = (\frac{2}{9}, 0, 0, 0, \frac{1}{9}, 0)^T.  $$   Since
$\|x^{(1)}\|_1 > \|x^{(3)}\|_1$ and $\|x^{(2)}\|_1 > \|x^{(3)}\|_1$,
by Theorem 3.2, the RSP of $A^T$ is impossible to hold at $x^{(1)}$
and $x^{(2)}. $ So we only need to check the RSP at $x^{(3)}. $
Taking $ y = (5, 5/3,0)^T$ yields $\eta= A^T y = (1, 1/3, -2/3,
-1/6, 1, -7/6)^T \in {\cal R} (A^T)$ where the first and fifth
components are 1, and all others are strictly less than 1. Thus the
RSP (\ref{RSP+}) holds at $x^{(3)},$ which (by Theorem 3.2) is the
unique least $\ell_1$-norm nonnegative solution to the linear
system.
   So  the
    $\ell_1$-method   solves the $\ell_0$-problem, although the $\ell_0$-problem  has
  multiple optimal solutions.

 The following corollary is an immediate consequence of Theorem 3.2,
which claims that when an $\ell_0$-problem has multiple sparsest
nonnegative optimal solutions, only one of them can satisfy the RSP
of $A^T.$

\vskip 0.08in

\textbf{Corollary 3.5 } \emph{For any underdetermined  system of
linear equations, there exists at most one sparsest nonnegative
solution satisfying the RSP (\ref{RSP+}).}

\vskip 0.08in

Example 3.4 and Theorem 3.2  show that $\ell_0$- and
$\ell_1$-problems can be equivalent provided that the RSP
(\ref{RSP+})
 is satisfied at an optimal solution to the $\ell_0$-problem, irrespective of the
multiplicity of   optimal solutions to the $\ell_0$-problem.  This analysis indicates that the success of the $\ell_1$-method
can happen not only for a subclass of linear systems in ${\cal
G}_1,$ but also for a wide range of linear systems in ${\cal G}_2.$
Since many existing conditions  imply the strong equivalence
 between
$\ell_0$- and $\ell_1$-problems, they can only explain the success
of $\ell_1$-methods when solving some $\ell_0$-problems in ${\cal
G}_1$. In other words, these strong equivalence  conditions cannot apply to the
$\ell_0$-problems in ${\cal G}_2$ which have multiple sparsest optimal
solutions, and hence they cannot
 interpret the  numerical success of the
$\ell_1$-method in these situations. However, the RSP-based analysis has shown that the  success of the
$\ell_1$-method may take  place
 not only  for those problems in  ${\cal
 G}_1$, but for a wide range of linear systems in ${\cal G}_2$ as well. So the  $\ell_1$-method  can solve a wider range of $\ell_0$-problems than what the strong equivalence theory can cope with. This does explain and clarify the gap between the  performance  of the
 $\ell_1$-method observed from numerical simulations and that indicated by existing strong equivalence conditions. Thus this analysis  provides  a certain  answer to the question (b) in Sect. 1.

\vskip 0.08in

\vskip 0.08in

\textbf{Remark 3.6 } It is worth noting
that our analysis method and results   can be easily generalized to
interpret the relationship between  $\ell_0$- and weighted
$\ell_1$-problems. More specifically, let us consider the
weighted $\ell_1$-problem \begin{equation} \label{W111}
\min\{\|Wx\|_1: Ax=b, x\geq 0\},\end{equation}  where
$W=\textrm{diag}(w)$ with $w\in R^n_+ $ and $w>0.$ By the nonsingular
linear transformation, $u=Wx, $ the above weighted $\ell_1$-problem
is equivalent to
\begin{equation} \label{uuuu} \min\{\|u\|_1: (AW^{-1})u=b, u\geq 0\}. \end{equation}  Clearly, $x$ is
the unique optimal solution to the  problem (\ref{W111}) if
and only if $u=Wx$ is the unique optimal solution to the
problem (\ref{uuuu}), and $u$ and $x$ have the same
supports. Thus any weighted $\ell_1$-problem with weight
$W=\textrm{diag}(w),$  where $w$ is a positive vector in $R^n_+,$  is
nothing but a standard $\ell_1$-problem with a scaled matrix
$AW^{-1}.$ As a result, applying Theorems 2.7
 to the problem (\ref{uuuu}),  we
 conclude that $u$ is the unique optimal solution to   (\ref{uuuu}) if and
 only if $ (AW^{-1})_{J_+(u)}$ has full column rank, and  there
 exists a vector $ \zeta \in {\cal R}((AW^{-1})^T)$ such that
 $\zeta_i=1$ for $u_i>0$ and $\zeta_i<1 $ for $u_i=0.$ By the one-to-one correspondence
 between the solutions of (\ref{W111}) and (\ref{uuuu}),
 and by transforming back to the weighted $\ell_1$-problem using  $u= Wx$ and $ \eta= W\zeta,$ we immediately conclude
 that \emph{$x$ is the unique optimal solution to the weighted
 $\ell_1$-problem (\ref{W111})
 if and only if (i) $A_{J_+}$  has full column rank where $J_+=\{i: x_i>0\},$ and (ii)
 there exists an $\eta \in {\cal R} (A^T)$ such that $\eta_i=w_i$
 for $x_i>0,$ and $\eta_i < w_i$ for $x_i=0.$}  We may call the above property (ii)
   the \emph{weighted RSP of $A^T$ at $x$}.  Thus the results developed in this paper
  can be easily generalized to the weighted $\ell_1$-method  for $\ell_0$-problems.

 \vskip 0.08in

\textbf{Remark 3.7}  The RSP-based analysis and results  can be also applied to the sparsest optimal solution of
the linear program (LP)
\begin{equation}\label{LLPP}   \min\{c^T x: Ax=b, x\geq 0\}.
\end{equation}
The sparsest optimal solution of (\ref{LLPP}) is meaningful. For
instance, in  production planning scenarios, the decision variables
$x_i\geq 0,$ $i=1, ..., n,$  represent what production
activities/events that should take place and how much resources
should be allocated to them in order to achieve an optimal objective
value. The sparsest optimal solution of a linear program provides
the smallest number of
  activities to achieve the optimal objective value. In many situations,
  reducing the number of activities is
vital for efficient planning, management and resource allocations.
We denote by $d^*$ the optimal value of (\ref{LLPP}), which can be
obtained by solving the LP by simplex methods, or interior point
methods. We assume that (\ref{LLPP}) is feasible and has a finite
optimal value $d^*.$  Thus the optimal solution set of the LP is
given by $\{x: ~ Ax=b,~ x\geq 0, ~c^Tx=d^*\}.$ So a sparsest optimal
solution to the LP is an optimal solution to the $\ell_0$-problem
\begin{equation} \label{LPl0} \min \left\{\|x\|_0: \left(
                     \begin{array}{c}
                       A \\
                       c^T \\
                     \end{array}
                   \right) x = \left(
                               \begin{array}{c}
                                 b \\
                                 d^* \\
                               \end{array}
                             \right), x\geq 0 \right\},
                             \end{equation}
associated with which is the $\ell_1$-problem
\begin{equation}\label{LPl1} \min \left\{\|x\|_1: \left(
                     \begin{array}{c}
                       A \\
                       c^T \\
                     \end{array}
                   \right) x = \left(
                               \begin{array}{c}
                                 b \\
                                 d^* \\
                               \end{array}
                             \right), ~ x\geq 0 \right\}.  \end{equation}
Therefore the developed results for sparsest nonnegative solutions
of linear systems in this paper
  can be directly applied to (\ref{LPl0})
 and (\ref{LPl1}). For instance,  from Theorems 2.7 and 3.2, we
 immediately have the following statements:  \emph{$x$ is the unique least $\ell_1$-norm
 optimal solution to the problem
   (\ref{LLPP}) if and only if    {\small $H =\left(
     \begin{array}{c}
       A_{J_{+}} \\
       c^T_{J_+}
     \end{array}
   \right)
$}  has full column rank, and  there  is a vector $\eta\in R^n$
satisfying
\begin{equation} \label{LPRSP+} \eta \in { \mathcal R} ( [A^T ,~ c]
), ~ \eta_i = 1  \textrm{ for all } i\in J_{+}, \textrm{ and }
\eta_i < 1 \textrm{ for all }i \notin J_+
\end{equation}
 where $J_{+} = \{ i: x_i >0\}. $   Moreover,   a sparsest optimal solution  to
  the problem  (\ref{LLPP}) is the unique least  $\ell_1$-norm   optimal solution to the LP
   if and only if the range space property (\ref{LPRSP+}) holds at this optimal solution. }
 Note  that a degenerated optimal solution has
been long studied since 1950s (see [47$-$48])   and the references
therein). It is well-known that finding a degenerated optimal
solution requires extra effort than nondegenerated ones. Finding the
most degenerated optimal solution or the sparsest optimal solution
becomes even harder. The RSP-based analysis provides a new understanding for the most
degenerated or the sparsest optimal solutions of LPs.

\section{Application to compressed sensing}

   One
of the tasks in compressed sensing is to exactly recover a sparse
vector (representing a signal or
 an image)
 via an underdetermined system of linear equations  [31, 33$-$35]. In this section, we
  consider the exact recovery of an unknown sparse nonnegative  vector
  $x^*$ by  $\ell_1$-minimization. For this purpose, we assume that an $m\times n $ $ (m<n)$ sensing
   matrix $A$ and the measurements $y=Ax^*$ are available.   A nonnegative solution $x $ of  the system $Ax=b$ is said to
  have a guaranteed recovery (or to be exactly recovered) by $\ell_1$-minimization   if $x$ is the
  unique least $\ell_1$-norm nonnegative solution to the linear system.
   To guarantee the success of recovery,  the current compressed
sensing theory
  assumes that the  matrix $A\in R^{m\times n} (m<n) $  satisfies
  some conditions  (e.g., the RIP or the NSP of order $2K$) which  imply the following properties:
  (i) $x^*$ is the  unique least $\ell_1$-norm nonnegative   solution to the
  system $Ax=y =Ax^*$ (where  the components of $y$ are
   measurements);
    (ii) $x^*$ is  the unique sparsest nonnegative   solution to
the system $Ax=y. $  So the underlying $ \ell_0$-  and $\ell_1$-problems
  must be strongly equivalent. Most of the recovering
conditions developed so far are for the so-called uniform recovery.

\subsection{Uniform recovery of sparse nonnegative vectors}
 The exact recovery of all $K$-sparse
nonnegative vectors (i.e., $\{x: ~ x\geq 0, ~\|x\|_0\leq K\}$) by a
single sensing matrix $A$ is called the \emph{uniform recovery} of
$K$-sparse nonnegative vectors. To develop a RSP-based recovery
theory, let us first introduce the following concept.

\vskip 0.08in

\textbf{Definition 4.1} (RSP of order $K$). \emph{Let $A$ be an
$m\times n$ matrix with $m<n. $  $A^T$ is said to satisfy the  range
space property of order $K$ if for any
  subset  $S\subseteq \{1,..., n\}$ with $|S|\leq
 K$,   $ {\mathcal R}(A^T) $ contains a vector $\eta$ such
 that $\eta_i =1 $ for all $ i \in S,$   and $ \eta_i  <1 $ for all $i\in S_c=\{1,2,..., n\}\backslash S.$}

 \vskip 0.08in
We first show that if $A^T$ has the RSP of order $K$, then $K$ must
be bounded by the  spark of $A$, denoted by Spark($A), $  which is
  the smallest number of columns of A that are linearly
dependent (see, e.g., [30, 38]).

\vskip 0.08in

\textbf{Lemma 4.2 }  \emph{If $A^T$ has the RSP of order $K$, then
any $K$ columns of $A$ are linearly independent, so $K< Spark(A).$}

\vskip 0.08in

\emph{Proof.} Let $S= \{s_1, ..., s_K\}, $ with $|S|=K ,$ be  an
arbitrary subset of $\{1,..., n\}. $  Suppose that $A^T$ has the RSP
of order $K.$  We now prove that $A_S$ has full column rank. It is
sufficient to show that $z_S=0$ is the only solution to $A_Sz_S=0. $
Indeed, let $A_Sz_S=0. $  Then $z=(z_S, z_{S_c} =0 ) \in R^n $ is in
the null space of $A. $   By the RSP of order $K$, there exists a
vector $\eta \in {\mathcal R}(A^T) $ such that every component of
$\eta_S$ is 1, i.e., $\eta_{s_i}=1$ for $i=1,..., K.$  By the
orthogonality of the null and range spaces, we have
\begin{equation} \label{zzss}    z_{s_1} + z_{s_2} + \cdots + z_{s_K} = z_S^T \eta_S =z^T \eta=0. \end{equation}
Now let $k$ be an arbitrary number with $1\leq k\leq K$, and
$S_k=\{s_1, s_2, ..., s_k\}\subseteq S. $ Since $|S_k|\leq |S|=K,$
it follows from  the definition of the RSP of order $K$, there exists a
vector $\widetilde{\eta} \in {\mathcal R}(A^T) $  with
$\widetilde{\eta}_{s_i}=1$ for every $i=1,..., k $ and
$\widetilde{\eta}_{j} <1 $ for every $j\notin S_k . $ By the
orthogonality again, it follows from $z^T\widetilde{\eta}=0$   that
$$ (z_{s_1}+ \cdots +z_{s_k}) + ( \widetilde{\eta}_{s_{k+1}}z_{s_{k+1}}
+ \cdots + \widetilde{\eta}_{s_K}z_{s_K}  )  =0 . $$ This is
equivalent to
$$ (z_{s_1}+ \cdots +z_{s_k}) +  (z_{s_{k+1}}+ \cdots +z_{s_K})   + [z_{s_{k+1}}(\widetilde{\eta}_{s_{k+1}}-1) +
\cdots + z_{s_K}  (\widetilde{\eta}_{s_K}-1)]  =0  $$ which,
together with (\ref{zzss}), implies that
$$ (\widetilde{\eta}_{s_{k+1}}-1) z_{s_{k+1}} + \cdots + (\widetilde{\eta}_{s_K}-1) z_{s_K}
 =0 $$ where $\widetilde{\eta}_{s_i}<1$ for $i= k+1, ..., K. $
Since such relations hold for every specified $k$ with $1\leq k\leq
K. $ In particular, for $k=K-1,$  the relation above is reduced to
$(\widetilde{\eta}_{s_K}-1) z_{s_K} =0 $ which implies that
$z_{s_K}=0 $ since $ \widetilde{\eta}_{s_K}<1.$  For $k=K-2$, the
relation above is of the form $$ (\widetilde{\eta}_{s_{K-1}}-1)
z_{s_{K-1}} + (\widetilde{\eta}_{s_K}-1) z_{s_K}  =0 $$ which,
together with $z_{s_K}=0$ and $ \widetilde{\eta}_{s_{K-1}} <1, $
implies that $ z_{s_{K-1}} =0.$ Continuing this process by
considering $k=K-3, ..., 1, $ we  deduce that all components of
$z_S$ are zero. Thus $A_S$ has full column rank. By the definition
of Spark$(A)$, we must have $K<\textrm{Spark} (A).$  \hfill $\Box$

 The  RSP of order $K$ can completely characterize  the uniform recovery  of all $K$-sparse nonnegative vectors by
 $\ell_1$-minimization, as shown by the next result.

\vskip 0.08in

\textbf{Theorem 4.3 } \emph{Let the measurements of the form $y= Ax$
be taken.  Then any $x\geq 0 $ with $\|x\|_0 \leq K$ can be exactly
recovered by the $\ell_1$-method (i.e., $\min\{\|z\|_1: Az=y, ~z\geq
0\}$) if and only if  $A^T$ has the RSP of order $K.$ }

\vskip 0.1in

\emph{Proof.}   Assume that the RSP of order $K$ is satisfied.  Let
$x^* \geq 0 $ be an arbitrary vector with $\|x^*\|_0 \leq K.$ Let
$S= J_+= \{i: x^*_i>0\}. $ Since $|S| =\|x^*\|_0 \leq K,$ by the RSP
of order $K$, there exists a vector $\eta \in {\mathcal R} (A^T) $
such that $\eta_i =1 $ for all $ i \in S,$ and $ \eta_i <1 $ for all
$i\in S_c.$ This implies that the  RSP (\ref{RSP+}) holds at
$x^*\geq 0.$ Moreover, it follows from Lemma 4.2 that $A_S$ has full
column rank. Hence, by Theorem 2.7,
    $x^* $ is  the unique least $\ell_1$-norm nonnegative solution to the system $Ax=y=Ax^*.$
    So $x^*$ can be exactly recovered by the $\ell_1$-method.

Conversely, assume that any $x\geq 0$ with $\|x\|_0\leq K$ can be
exactly recovered by  the $\ell_1$-method. We now prove that the RSP
of order $K$ must be satisfied. Let  $S=J_+=\{i: x_i>0\}.$ Under the
assumption, $x$ is  the unique optimal solution to the
$\ell_1$-problem $$\min\{\|z\|_1: Az= y=Ax, z\geq 0\}.$$  By Theorem
2.7, the RSP (\ref{RSP+}) holds at $x,$  i.e., there exists a vector
$\eta\in {\cal R}(A^T)$ such that $\eta_i =1 $ for all $i\in S=J_+
$, and $\eta_i< 1$ otherwise. Since $x$ can be any $K$-sparse
nonnegative vectors, this implies that $S=J_+$ can be any subset  of
$  \{1, ..., n\}$ with $|S|\leq K, $ and for every such a subset
there exists accordingly a vector $\eta$ satisfying the above
property. By Definition 4.1, $A^T$ has the RSP of order $K.$  \hfill
$\Box$

  Let  $a_j,$ $ 1\leq j\leq n,$ be the columns of $A $ and let
$a_0=0.$ Let $P$ denote the convex hull of $a_j, $ $ 0\leq j\leq n.$
Donoho and Tanner [4] introduced the following concept: The polytope
$P$ is \emph{outwardly $K$-neighborly} if every subset of $K$
vertices not including $a_0=0$ spans a face of this polytope.  They
have shown
 that the polytope $P$ is outwardly $K$-neighborly if and only if
 any nonnegative solution $x$ to the system $Ax=b$
  with $\|x\|_0\leq K$ is the unique optimal solution to the
  $\ell_1$-problem. In other words, the outwardly $K$-neighborly
  property is a full geometric characterization of the uniform recovery of all $K$-sparse nonnegative vectors.
Recently, Khajehnejad et al [6] characterized the uniform recovery
by using the property of ${\cal N}(A),$ the
  null space  $A.$  They have showed that all nonnegative $K$-sparse
  vector can be exactly recovered if and only if
 for every vector $w\not=0$ in   ${\cal N}( A)$, and every
 index set $S\subseteq \{1, ..., n\}$ with $|S|=K$ such that
 $w_{S_c} \geq 0,$ it holds that $e^T w >0.$ As shown by Theorem 4.3,
  the RSP of order $K$ introduced in this section provides an alternative full characterization of the uniform
recovery of all $K$-sparse vectors.  Clearly,  from different perspectives, all the above-mentioned  properties
 (outwardly $K$-neighborly, null space, and range space of order $K$)
  are equivalent since  each of these properties is a necessary and sufficient condition for (i.e., equivalent to) the uniform
  recovery of all $K$-sparse vectors.  As a result, if a matrix satisfies one of these conditions,
  it  also  satisfies  the other two conditions.  However,  directly checking
  these conditions is difficult. Checking the outwardly
  $k$-property needs to check all possible $K$ out of $n$ vertices
  and verify whether all such subsets can span  a face of the polytope
  defined by $A.$ To check the null space property, we need to
  to check all possible vectors in the null space of $A,$ and for every such a vector, we need to verify  the required property
for all possible $K$ components of the vector.
 Similarly, for the range space property, we need to
  verify  the individual RSP holds at
 every  subset $S\subset \{1,2,..., n\} $  with cardinality
 $ |S| \leq K .$  Clearly, none of these properties is easier to check than the others
 for general instances of matrices.


We now close this section by stressing the difference between the
RSP of order $K$ and the RSP (\ref{RSP+}). Such a difference can be
easily seen from the following result.

\vskip 0.08in

 \textbf{Corollary 4.4 }  \emph{If $A^T$ has the RSP of order
$K$, then any $\widehat{x}\geq 0$ with $\|\widehat{x}\|_0 \leq K$ is
both the unique least $\ell_1$-norm nonnegative solution and the
unique sparsest nonnegative solution to the linear system
$Ax=y=A\widehat{x}.$}

 \vskip 0.08in

 \emph{Proof.}  By Theorem 4.3,
under the RSP  of order $K,$ any $\widehat{x} \geq 0$ with
$\|\widehat{x}\|_0\leq K $ can be exactly recovered by
$\ell_1$-minimization, i.e, $\widehat{x}$ is the unique least
$\ell_1$-norm nonnegative solution to the   system
$Ax=y=A\widehat{x}.$ We now prove that $\widehat{x}$ is also the
sparsest nonnegative solution to this system. Assume that there
exists another solution $z\geq 0 $ such that $\|z\|_0\leq
\|\widehat{x}\|_0.$ Let $S=\{i: z_i>0\}.$ Since $|S| =\|z\|_0 \leq
\|\widehat{x}\|_0\leq K,$ by the RSP of order $K$, there exists an
$\eta\in {\cal R} (A^T) $ such that $\eta_i =1 $ for all $i\in S$,
and $\eta_i<1 $ for all $i\in S_c. $ Thus the individual RSP
(\ref{RSP+}) holds at $z. $ By Lemma 4.2, any $K$ columns of $A$ are
linearly independent. Since the number of the columns of $A_S$,
where $S=\{i: z_i>0\}$, is less than $K$, this implies that $A_S$
has full column rank. By Theorem 2.7, $z$ is also the unique least
$\ell_1$-norm nonnegative solution to the system
$Ax=y=A\widehat{x}.$ Thus $z=\widehat{x},$ which implies that
$\widehat{x}$  is   the unique sparsest nonnegative solution to this
system. \hfill $\Box$

This result shows that the RSP of order $K$ is  more
restrictive than the individual RSP (\ref{RSP+}) which is defined at
a single point. The former requires that the RSP (\ref{RSP+}) hold
at every $K$-sparse nonnegative solution. By contrast,  the
individual RSP (\ref{RSP+}) is only a local property, and it does
not imply that the underlying linear system   has a unique sparsest
nonnegative solution, as we have shown in Sect. 3.

\subsection{Non-uniform recovery of sparse nonnegative vectors}

The purpose of the uniform recovery is to  recover all $k$-sparse
vectors. So some strong assumptions (such as the RIP, NSP and RSP of
certain orders) must be imposed on the matrix. These strong
assumptions   imply that the unknown
sparse vector $x$ must be the unique optimal solution to both
$\ell_0$- and $\ell_1$-problems (hence, the strong equivalence
between these two problems are required by the uniform recovery). In
this subsection, we extend the uniform recovery theory to
non-uniform ones by using the RSP. The non-uniform recovery of sparse
signals has been investigated by some researchers (see e.g., [4, 49]). From a geometric perspective, Donoho
and Tanner [4] introduced the so-called weak neighborliness
conditions for the non-uniform recovery by
 $\ell_1$-minimization, and they have shown under such a condition that most
nonnegative $K$-sparse vectors can be exactly recovered by the
$\ell_1$-method. Ayaz and Rauhut [49] focused on the non-uniform
recovery of signals with given sparsity and given signal length by
$\ell_1$-minimization,  and they have provided the number of samples
required to recover such signals with gaussian and subgaussian
random matrices. In what follows, we introduce the
so-called weak RSP of order $K, $  a range
space property of $A^T,$  which provides a non-uniform recovery condition for
some vectors which may have high sparsity level, going beyond the
scope of normal uniform recoveries.

  Given a sensing matrix $A$,  Theorem 2.7
claims that  a vector $x^*$  can be exactly recovered by
$\ell_1$-minimization  provided that  the RSP(\ref{RSP+}) hold at
$x^*$ and that the matrix $A_{J_+}$, where $J_+ =\{i: x^*_i>0\},$
 has full-column rank. Such an $x^*$ is not necessarily  the
unique sparsest nonnegative solution to the linear system
 as shown by Example 3.4, and it may not even be a
sparsest nonnegative solution as well. For instance, let
   $$A = \left(
  \begin{array}{crrrr}
    6 & 4 & 1.5 & 4 & -1 \\
    6 & 4 & -0.5 & 4 & 0 \\
    0 & -2 & 31.5 & -1 & -1.5 \\
  \end{array}
\right), ~   y = \left(
                   \begin{array}{r}
                     4 \\
                     4 \\
                     -1 \\
                   \end{array}
                 \right) = Ax^*
$$ where $x^*= (1/3, 1/2,
0, 0, 0)^T. $  It is easy to see that  $\widetilde{x}=(0, 0, 0, 1,
0)^T$ is the unique sparsest nonnegative solution to the system
$Ax=y, $ while $x^*$ is the unique least $\ell_1$-norm nonnegative
solution to the system $Ax=y. $ Although $x^* $ is not the sparsest
nonnegative solution, it can be exactly recovered by the
$\ell_1$-method. Because of this, it is interesting to  develop a
recovery theory without requiring that the targeted unknown sparse
vector be a sparsest or be the unique sparsest solution to a linear
system. This is also motivated by some practical applications. In
fact, a real sparse signal   may not be  sparse enough to be
recovered by the uniform recovery, and  partial information for the
unknown sparse vector may be  available in some situations, for
example, the support of an unknown vector may be known. The concept
of the RSP of  order $K$ can be  adapted to handle these cases. So
we introduce the following concept.

\vskip 0.08in

\textbf{Definition 4.5}(WRSP of order $K$)~   \emph{Let $A$ be an
$m\times n$ matrix with $m<n.$   $A^T$ is said to satisfy the weak
range space property  (WRSP) of order $K$ if the following two
properties are satisfied:}

  (i) \emph{There exists a  subset $S\subseteq \{1,..., n\}$
such that $ |S|=
 K $ and  $ A_{S} $ has full column rank.}

(ii) \emph{For any subset $S \subseteq \{1,..., n\}$ such that $ |S|\leq
 K $ and $ A_{S} $ has full column rank,  the space $ {\mathcal R}(A^T) $
contains a vector $\eta$ such
 that $\eta_i =1 $ for  $ i \in S,$ and $ \eta_i <1
 $ otherwise.}

\vskip 0.08in

 The WRSP of order $K$ only requires that the
individual RSP hold for those subsets $S\subseteq \{1,..., n\}$ with
$|S|\leq K$ and that $
      A_{S}$  be full-column-rank,  while the RSP of order $K$ requires
that the individual RSP
  hold for any subset $S\subseteq \{1,...,
n\}$ with $|S|\leq K. $   So the WRSP of order $K$ is less restrictive than
the RSP of order $K.$   By Theorem 2.6, we have the following
result.

\vskip 0.08in

 \textbf{Theorem 4.6 }  ~\emph{Let the measurements of the form $y= Ax$
be taken. Suppose that there exists a subset $S\subseteq \{1,...,
n\}$ such that $ |S|=
 K $ and $
      A_{S}$ has full column rank.  Then $A^T$ has the WRSP of order $K$ if and only if
any $x\geq 0,$ satisfying that $\|x\|_0 \leq K$ and $
      A_{J_+}  $ has full-column-rank where $J_+=\{i: x_i>0\},$  can be exactly
recovered by the $\ell_1$-minimization  $ \min \{\|z\|_1: Az=y =Ax ,
z\geq 0\}. $ }

\vskip 0.08in

\emph{Proof.}  Assume that $A^T$ has the WRSP of order $K.$ Let $x$
be an arbitrary nonnegative vector such that $\|x\|_0\leq K$ and $
      A_{J_+}  $ has full-column-rank, and let $S=J_+=\{i: x_i>0\}.$   Since $A^T$ has the WRSP of order $K,$
  there exists an $\eta\in {\mathcal R}(A^T) $ such
that $\eta_i= 1$ for $i\in S=J_+,$ and $\eta_i<1$ otherwise. This
implies that  the RSP(\ref{RSP+}) holds at $x.$ Since $
      A_{J_+} $ has full column rank, by
Theorem 2.7, $x$ must be the unique least $\ell_1$-norm nonnegative
solution to the  system $Az=y $ $ (=Ax).$ In other words, $x$
can be  recovered by $\ell_1$-minimization.   Conversely, we
assume that any $x\geq 0,$ satisfying that $\|x\|_0 \leq K$ and $
A_{J_+} $ has full-column-rank,  can be exactly recovered by
$\ell_1$-minimization. We now prove that $A^T$ must have the WRSP of
order $K.$ In fact, let $x\geq 0$ be a vector such  that $\|x\|_0
\leq K$ and $
      A_{J_+}$  has full-column-rank.
  Denote by $S=J_+=\{i: x_i>0\}.$ Since $x$ can be  recovered
by the $\ell_1$-method, it is the unique least $\ell_1$-norm
nonnegative solution to the system $Az=y=Ax.$ By Theorem 2.7, the
RSP (\ref{RSP+}) holds at $x$, i.e., there exists an $\eta\in
{\mathcal R} (A^T) $ such that $\eta_i=1$ for $i\in J_+ =S,$ and
$\eta_i<1$ otherwise. Since $x$ can be any vector such that
$\|x\|_0\leq K $ and $
      A_{J_+} $  has full column rank, this implies that the condition
(ii) of Definition 4.5 holds, thus $A^T$ has the WRSP of order $K.$
 \hfill $\Box$

   We may further relax the concept
of the RSP and WRSP, especially when
  partial information  available to the unknown vector.  For instance, when $\|x\|_0 = K $ is known,
   we may introduce the next two concepts.

\vskip 0.08in

\textbf{Definition 4.7  } ~  (PRSP of order $K$).  \emph{We say that
$A^T$ has the partial range space property (PRSP) of order $K$ if
for any subset $S$ of $\{1,..., n\}$ with $ |S|=
 K$, the range space $ {\mathcal R}(A^T) $ contains a vector $\eta$ such
 that $\eta_i =1 $ for all $ i \in S,$ and $ \eta_i <1 $ otherwise. }

\vskip 0.08in

 \textbf{Definition 4.8 } ~ (PWRSP of order $K$).
\emph{$A^T$ is said to have partial weak range space property
(PWRSP) of order $K$ if
 for any  subset $S\subseteq \{1,..., n\}$ such
that $ |S|=
 K $ and  $
      A_{S}$  has full column rank,    $ {\mathcal R}(A^T) $ contains a vector $\eta$ such
 that $\eta_i =1 $ for all $ i \in S,$  and $ \eta_i <1 $
 otherwise.
}

\vskip 0.08in

Different from the RSP of order $K$, the PRSP of order $K$ only
requires that  the individual RSP hold  for the subset $S$ with
$|S|=K. $  Similarly, the PWRSP of order $K$ is also less
restrictive than WRSP. Based on such definitions, we have the next
result which follows from Theorem 2.7 straightaway.

\vskip 0.08in

\textbf{Theorem 4.9 } ~ (i)   \emph{The matrix $A^T$ has the partial
range space property (PRSP) of order $K$ if and only if any $x\geq
0,$ with $\|x\|_0 = K,$ can be exactly recovered by the
$\ell_1$-minimization  $ \min \{\|z\|_1: Az=y = Ax, z\geq 0\}. $}

(ii) \emph{$A^T$ has the PWRSP of order $K $ if and only if any $x
\geq 0,$ satisfying that  $\|x\|_0 =K$ and $
      A_{J_+}  $  has full-column-rank where $J_+= \{i: x_i>0\},$  can be exactly
recovered by the $\ell_1$-minimization  $ \min \{\|z\|_1: Az=y = Ax
, z\geq 0\}. $}

\vskip 0.08in

 When $ A_{S} $ has full column rank, we have $|S|\leq m.$ Thus
     the   WRSP and PWRSP of order
   $K$ imply that $K\leq m.$  Moreover, the   PRSP of order $K$
   implies that  $K<\textrm{Spark}(A). $ In fact, the proof of  this fact is identical to that of Lemma
   4.1.
 Theorems 4.6 and 4.9(ii) indicate that a portion of vectors with $\|x\|_0 \leq m $ can be
 recovered if a sensing matrix satisfies certain properties
milder than the RSP of order $K$ (and thus milder than the RIP and the NSP
of order $2K$). Since the PRSP, WRSP and PWRSP of order $K$ do not
require that the individual RSP hold for all subsets $S$ with
$|S|\leq K,$ by Theorem 4.3, these properties are non-uniform
recovering conditions developed through certain range space properties of
$A^T.$

 \section{Conclusions}

 Through the range space property, we
have characterized the conditions for  an $\ell_1$-problem   to
have a unique optimal solution,   and  for $K$-sparse vectors to be
uniformly or  non-uniformly recovered by the $\ell_1$-method.  We have shown the following
 main results: (i) A nonnegative vector  is the
unique optimal solution to the $\ell_1$-problem if and only if the
RSP holds at this vector, and the associated submatrix $A_{J_+}$ has
full column rank;  (ii) All $K$-sparse vectors can
be exactly recovered by a single sensing matrix    if and only if the transpose of this matrix
has the RSP of order $K;$ (iii) $\ell_0$- and $\ell_1$-problems are
equivalent if and only if the   RSP  holds at an optimal
solution of the $\ell_0$-problem. The
  RSP originates naturally from the strict complementarity
property of linear programming problems.     The RSP-based  analysis has indicated that the
uniqueness of optimal solutions of $\ell_0$-problems is not the
reason for the problem being computationally tractable, and the
multiplicity of optimal solutions of $\ell_0$-problems is also
not the reason for the problem being hard. The RSP may hold in both
situations, and hence an $\ell_0$-problem can be successfully solved by the $\ell_1$-method  in both situations, provided that the RSP is satisfied at a solution of the  $\ell_0$-problem.
Thus the relationship between $\ell_0$- and $\ell_1$-problems and the numerical performance of the $\ell_1$-method can be broadly interpreted via the RSP-based analysis.\\

\noindent
  \textbf{References}
\begin{enumerate}

\item[1.] Chen, D., Plemmons, R.: Nonnegativity constraints in
numerical analysis.  Symposium on the Birth of Numerical Analysis
(2007)

\item[2.] Bardsley, J.,  Nagy, J.: Covariance-preconditioned iterative
methods for nonnegativity constrained astronomical imaging.  SIAM J.
Matrix Anal. Appl. \textbf{27}, 1184$-$1198 (2006)

\item[3.] Bruckstein, A.M., Elad, M.,  Zibulevsky, M.: On the
uniqueness of nonnegative sparse solutions to underdetermined
systems of equations. IEEE Trans. Inform. Theory \textbf{54},
4813$-$4820 (2008)

\item[4.] Donoho, D.L.,  Tanner, J.: Sparse nonnegative
solutions of underdetermined linear equations by linear programming.
Proc.  Natl. Acad. Sci. \textbf{102}, 9446$-$9451 (2005)

\item[5.] Donoho, D.L., Tanner, J.:  Counting the faces of
randomly-projected hypercubes and orthants with applications.
Discrete Comput. Geom. \textbf{43}, 522$-$541 (2010)

\item[6.] Khajehnejad, M.A.,  Dimakis, A.G.,  Xu, W.,
Hassibi, B.:  Sparse recovery of nonnegative signals with minima
expansion.  IEEE Trans. Signal Proc. \textbf{59}, 196$-$208 (2011)

\item[7.]  O'Grady, P.D.,   Rickard, S.T.: Recovery of
nonnegative signals from compressively samples observations via
nonnegative quadratic programming. SPASR'09 (2009)

\item[8.] Slawski, M.,  Hein, M.: Sparse recovery by
thresholded nonnegative least squares. Technical Report, Saarland
University, Germany (2011)

\item[9.] Vo, N.,  Moran, B., Challa, S.:  Nonnegative-least-square
classifier for face recognition. Proc. of the 6th Inter. Symposium
on Neural Networks: Advances in Neural Networks, 449-456 (2009)

\item[10.] Wang, M.,  Xu, W., Tang,  A.: A unique nonnegative
solution to an underdetermined system: from vectors to matrices.
IEEE Trans. Signal Proc. \textbf{59}, 1107$-$1016 (2011)

\item[11.] Bradley, P.S.,  Mangasarian, O.L., Rosen, J.B.:
Parsimonious least norm approximation. Comput. Optim. Appl.
\textbf{11}, 5$-$21 (1998)

\item[12.]  Bradley, P.S.,  Fayyad, U.M., Mangasarian, O.L.:
Mathematical programming for data mining: formulations and
challenges. INFORMS J. Computing \textbf{11}, 217$-$238 (1999)

\item[13.] He, R.,  Zheng, W., Hu, B., Kong, X.:  Nonnegative
sparse coding for discriminative semi-supervised learning. Proc. of
IEEE conference on CVPR, 2849$-$2856 (2011)

\item[14.] Mangasarian, O.L.:  Minimum-support solutions of
polyhedral concave programs. Optimization \textbf{45}, 149$-$162
(1999)

\item[15.] Mangasarian, O.L., Machine learning via polydedral
concave minimization.  Applied Mathematics and Parallel
Computing-Festschrift for Klaus Ritter (H. Fischer, B. Riedmueller
and S. Schaeffler eds.),  Springer, Heidelberg, 175$-$188  (1996)

\item[16.]  Szlam, A.,  Guo, Z.,  Osher, S.:  A split Bregman
method for nonnegative sparsity penalized least squares with
applications to hyperspectral demixing. IEEE Inter. Conf. on Imaging
Processing, 1917$-$1920 (2010)

\item[17.]  Slawski, M., Hein,M.:  Sparse recovery for protein
mass spectrometry data.  Technical Report, Saarland University,
Germany (2010)

\item[18.] Hoyer, P.O.:  Nonnegative sparse coding.  Proc. of the
12th IEEE Workshop, 557$-$565 (2002)

\item[19.] Wang, F., Li, P.: Compressed nonnegative sparse
coding.  IEEE 10th Int. Conf. on Data Mining (ICDM), 1103$-$1108
(2010)

\item[20.]  Recht, B.,   Fazel, M.,   Parrilo, P.A.:  Guaranteed
minimum-rank solutions of linear matrix equations via nuclear norm
minimization.  SIAM Rev. \textbf{52},  471$-$501 (2010)

\item[21.] Zhao, Y.B.:  An approximation theory of matrix rank
minimization and its application to quadratic equations. Linear
Algebra  Appl. \textbf{437}, 77$-$93 (2012)

\item[22.] Lee, D., Seung, H.:  Learning the parts of objects
by nonnegative matrix factorization.  Nature \textbf{401}, 788$-$791
(1999)

\item[23.] Park, H.,  Kim, H.: Nonnegative matrix factorization
based on alternating nonnegativity constrained least squares and the
active set method. SIAM J. Matrix Anal. Appl. \textbf{30}, 713$-$730
(2008)

\item[24.]  Pauca, P., Piper, J., Plemmons, R.:  Nonnegative matrix
factorization for spectral data analysis.  Linear Algebra Appl.
\textbf{416}, 29-47 (2006)

\item[25.]  Natarajan, B.K.:  Sparse approximate solutions to linear
systems. SIAM J. Comput. \textbf{24}, 227$-$234 (1995)

\item[26.]  Cand$\grave{\textrm{e}}$s, E.,  Wakin, M.B.,
Boyd,  S.P.: Enhancing sparsity by reweighted $\ell_1$ minimization.
J. Fourier Anal. Appl. \textbf{14},  877$-$905 (2008)

\item[27.]  Zhao, Y.B.,  Li, D.:  Reweighted $\ell_1$-minimization for sparse solutions to
underdetermined linear systems. SIAM J. Optim.  \textbf{22},
1065$-$1088 (2012)

\item[28.]  Wen, Z.,  Yin, W., Goldfarb, D.,  Zhang, Y.: A fast
algorithm for sparse reconstruction based on shrinkage, subspace
optimization, and continuation. SIAM J. Sci. Comput. \textbf{32},
1832$-$1857 (2010)

\item[29.] Yin, W.,   Osher, S.,  Goldfarb, D.,  Darbon, J.:
Bregman iterative algorithms for $\ell_1$ minimization with
applications to compressed sensing. SIAM J. Imaging Sci. \textbf{1},
143$-$168 (2008)

\item[30.] Bruckstein, A.M., Donoho, D.L.,  Elad, M.:  From sparse
solutions of systems of equations to sparse modeling of signals and
images. SIAM Rev. \textbf{51}, 34$-$81 (2009)

\item[31.]  Cand$\grave{\textrm{e}}$s,  E.,  \emph{Compressive sampling}. Proceedings of the
International Congress of Mathematicians, Madrid, Spain (2006)

\item[32.]  Cand\`es, E.,   Romberg, J.,  Tao, T.:  Robust
uncertainty principles: Exact reconstruction from highly incomplete
frequency information. IEEE Trans. Inform. Theory \textbf{52},
489$-$509 (2006)

\item[33.]   Donoho, D. L.:  Compressed sensing. IEEE Trans. Inform.
Theory \textbf{52}, 1289$-$1306 (2006)

 \item[34.]  Elad, M.:  Sparse and Redundant Representations: From Theory to Applications in Signal
and Image Processing.  Springer, New York (2010)

\item[35.]  Cand\`es, E.: The restricted isometry property and its
implications for compressed sensing.   Comptes Rendus de
l'Acad\'emie des Sciences, Paris, S\'erie I, \textbf{346}, 589$-$592
(2008)

\item[36.]  Cand$\grave{\textrm{e}}$s, E.,  Tao, T.: Decoding by linear
programming. IEEE Trans. Inform. Theory \textbf{51}, 4203$-$4215
(2005)

\item[37.]  Daubechies, I.,  DeVore, R.,  Fornasier, M.,
G\"unt\"urk, S.:  Iteratively reweighted squares minimization for
sparse recovery. Commun. Pure. Appl. Math. \textbf{63}, 1-38 (2010)

\item[38.]  Donoho, D.L.,  Elad, M.,  Optimality sparse
representation in general (non-orthogonal) dictionaries via $\ell_1$
minimization. Proc.  Natl. Acad. Sci. \textbf{100}, 2197$-$2202
(2003)

\item[39.]  Elad, M., Bruckstein, A.: A generalized uncertainty
principle and sparse representation in pairs of bases. IEEE Trans.
Inform. Theory \textbf{48}, 2558$-$2567 (2002)

 \item[40.]  Feuer, A.,   Nemirovsky, A.:  On sparse representation
in pairs of bases.  IEEE Trans. Inform. Theory  \textbf{49},
1579$-$1581 (2003)

\item[41.]  Fuchs, J.J.: On sparse representations in arbitrary
redundant bases. IEEE Trans. Inform. Theory \textbf{50}, 1341$-$1344
(2004)

\item[42.]  Tropp, J.A.:   Greed is good: Algorithmic results for sparse approximation.
 IEEE Trans. Inform. Theory \textbf{50}, 2231$-$2242 (2004)

\item[43.]  Juditski, A.,  Nemirovski, A.:  On verifiable
sufficient conditions for sparse signal recovery via $\ell_1$
minimization. Math. Program., Ser. B \textbf{127}, 57$-$88 (2011)

 \item[44.] Juditsky, A.,  Kilinc Karzan, F., Nemirovski, A.: Verifiable
conditions of $\ell_1$-recovery for sparse signals with sign
restrictions. Math. Program. Ser. B \textbf{127}, 89$-$122 (2011)

 \item[45.]  Gribonval, R.,  Nielsen, M.:  Sparse decompositions in
unions of bases. IEEE Trans. Inform. Theory \textbf{49}, 3320$-$3325
(2003)

 \item[46.]  Schrijver, A.: Theory of Linear and Integer
Programming.  Wiley, New York (1986)

 \item[47.]  Greenberg, H.J.:  An analysis of degeneracy.  Naval
Res. Logist. Quart. \textbf{33}, 635$-$655 (1986)

 \item[48.]  Greenberg, H.J.: The use of the optimal partition in a
linear programming solution for postoptimal analysis.  Oper. Res.
Lett. \textbf{15},  179$-$185 (1994)

\item[49.]  Ayaz, U.,   Rauhut, H.:  Nonuniform sparse recovery with
subgaussian matrices.  Technical Report, University of Bonn, Germany
(2011)

\end{enumerate}

 \end{document}